\documentclass[12pt]{article}
\usepackage{amsmath}
\usepackage{graphicx}
\usepackage{amscd}
\usepackage{amsfonts}
\usepackage{amssymb}
\newtheorem{theorem}{Theorem}[section]
\newtheorem{proposition}[theorem]{Proposition}
\newtheorem{lemma}[theorem]{Lemma}
\newtheorem{corollary}[theorem]{Corollary}

\begin{document}

\title{Remarks on the general Funk-Radon transform and thermoacoustic tomography}
\author{V. P. Palamodov}
\date{07.01.2007}
\maketitle

\section{Introduction}

The generalized Funk transform is a integral transform acting from densities
on a manifold $X$ to functions defined on a family $\Sigma$ of hypersurfaces
in $X.$ The dual operator $M%
{{}^\circ}%
$ is again a Funk transform which is defined for densities on the manifold
$\Sigma.$ We state two-side estimates for the operator $M$ and a description
of the range of the Funk transform operator $M$ and approximation theorem for
the kernel of the dual operator. This operator is similar to the integral
transform related to a double fibration in the sense of Guillemin \cite{G} and
some results can be extracted from his theory. Other results can be
generalized for the general double fibration.

\section{Geometry}

Let $X$ and $\Sigma$ be smooth $n$-manifolds, $n>1$ and $F$ be a smooth closed
hypersurface in $X\times\Sigma$ such that

(\textbf{i}) the projections $p:F\rightarrow X$ and $\pi:F\rightarrow\Sigma$
have rank $n.$ This condition implies that the sets $F\left(  \sigma\right)
=\pi^{-1}\left(  \sigma\right)  ,\sigma\in\Sigma$ and $F\left(  x\right)
=p^{-1}\left(  x\right)  ,x\in X$ are hypersurfaces in $X,$ respectively, in
$\Sigma.$ We call $F$ \textit{incidence} manifold. It can be defined locally
by the equation $I\left(  x,\sigma\right)  =0,$ where $I$ is a smooth function
in $X\times\Sigma$ such that $\mathrm{d}I\neq0$. If the hypersurface $F$ is
cooriented in $X\times\Sigma$, a function $I$ can be chosen globally; we call
it \textit{incidence} function. If $F$ is not cooriented, one can choose local
incidence functions $I_{\alpha},\alpha\in A$ such that $I_{\beta}=\pm
I_{\alpha}$ in the domain, where both functions $I_{\alpha}$, $I_{\beta}$ are
defined. The additional condition is

(\textbf{ii}) the mapping $q:F\rightarrow G^{n-1}\left(  X\right)  $ is a
local diffeomorphism, where $G^{n-1}\left(  X\right)  =\cup_{X}G_{x}^{n-1}$
and $G_{x}^{n-1}$ means the manifold of $n-1$-subspaces in the tangent space
$T_{x}$ of $X$ at $x;$ $q\left(  x,\sigma\right)  \doteq\left(  x,H\right)  ,$
where $H$ denotes the tangent hyperplane to $F\left(  \sigma\right)  $ at $x$.
It follows that for any point $x\in X$ and any tangent hyperplane $H\subset
T_{x}$ there is locally only one hypersurface $F\left(  \sigma\right)  $ that
contains $x$ and is tangent to $H.$

\begin{proposition}
\label{p1}The conditions (\textbf{i-ii}) are equivalent to the inequality
$\mathrm{\det\,}\Phi\neq0$ in $F,$ where
\[
\Phi=\left(
\begin{array}
[c]{cccc}%
\frac{\partial^{2}I}{\partial x_{1}\partial\sigma_{1}} & ... & \frac
{\partial^{2}I}{\partial x_{1}\partial\sigma_{n}} & \frac{\partial I}{\partial
x_{1}}\\
&  &  & \\
... & ... & ... & ...\\
&  &  & \\
\frac{\partial^{2}I}{\partial x_{n}\partial\sigma_{1}} & ... & \frac
{\partial^{2}I}{\partial x_{n}\partial\sigma_{n}} & \frac{\partial I}{\partial
x_{n}}\\
&  &  & \\
\frac{\partial I}{\partial\sigma_{1}} & ... & \frac{\partial I}{\partial
\sigma_{n}} & I
\end{array}
\right)  ,
\]
where $x_{1},...,x_{n}$ and $\sigma_{1},...,\sigma_{n}$ are local coordinates
in $X$ and in $\Sigma,$ respectively.
\end{proposition}

\textsf{Proof.} Suppose that $\det\Phi\neq0.$ Then also $\mathrm{d}_{x}I\neq0
$ and $\mathrm{d}_{\sigma}I\neq0$ which implies (\textbf{i}). Choose a point
$\left(  x_{0},\sigma_{0}\right)  $ and take a tangent vector $\theta$ to
$F\left(  \sigma_{0}\right)  $ at $x_{0}$. The vector $\left(  \theta
,0\right)  $ is tangent to $F$ at $\left(  x_{0},\sigma_{0}\right)  $ and the
map $q$ is well defined. Change the coordinates $x_{i}$ and $\sigma
_{i},i=1,...,n$ in such a way that $\partial I\left(  x_{0},\sigma_{0}\right)
/\partial\sigma_{i}=\partial I\left(  x_{0},\sigma_{0}\right)  /\partial
x_{i}=0,i=2,...,n$ at this point. We have then
\[
\det\Phi=\frac{\partial I}{\partial x_{1}}\frac{\partial I}{\partial\sigma
_{1}}\mathrm{\det}\Psi,\;\Psi\doteq\left\{  \frac{\partial^{2}I}{\partial
x_{i}\partial\sigma_{j}}\right\}  _{i,j=2}^{n}.
\]
The inequality $\mathrm{\det}\Psi\neq0$ implies that the forms $\partial
\mathrm{d}_{x}I\left(  x_{0},\sigma\right)  /\partial\sigma_{2},...,\partial
\mathrm{d}_{x}I\left(  x_{0},\sigma\right)  /\partial\sigma_{n}$ are
independent. This means that the fields $\partial/\partial\sigma
_{2},...,\partial/\partial\sigma_{n}$ do not move the point $x_{0}$ but rotate
the tangent hyperplane to $F\left(  \sigma_{0}\right)  $ at $x_{0}$ whereas
the field $\partial$/$\partial\sigma_{1}$ move the point $x_{0}.$ This yields
(\textbf{ii}). The inverse statement can be proved on the same lines.
$\blacktriangleright$

It follows that the properties (\textbf{i-ii}) are symmetric with respect to
$X$ and $\Sigma$. The next condition is not symmetric:

(\textbf{iii}) the projection $p:F\rightarrow X$ is proper. If $F$ satisfies
(\textbf{ii}) and (\textbf{iii})\textbf{ }and the hypersurface $F\left(
x\right)  $ is not empty for a point $x\in X,$ then the mapping $q:F\left(
x\right)  \rightarrow G_{x}^{n-1}$ is surjective, since the manifold
$G_{x}^{n-1}$ is connected.

\section{The Funk transform}

Consider a hypersurface $F$ as above defined by a incidence function $I$ that
fulfils (\textbf{i}). Define the Funk (or Minkowski-Funk-Radon) transform for
densities $\mathrm{f}$ in $X$ with compact support by means of the integral
\begin{equation}
M\mathrm{f}\left(  \sigma\right)  \doteq\lim_{\varepsilon\rightarrow0}\frac
{1}{2\varepsilon}\int_{\left|  I\left(  \cdot,\sigma\right)  \right|
\leq\varepsilon}\mathrm{f}=\int_{F\left(  \sigma\right)  }\frac{\mathrm{f}%
}{\mathrm{d}_{x}I},\;\sigma\in\Sigma, \label{1}%
\end{equation}
where $\mathrm{\omega}=\mathrm{f}/\mathrm{d}_{x}I$\ is a $n-1$-form such that
$\mathrm{d}_{x}I\wedge\mathrm{\omega}=\mathrm{f.}$ This form is defined up to
a term $\mathrm{d}_{x}I\wedge\mathrm{\chi},$ where $\mathrm{\chi}$ is a
$n-2$-form. Therefore the restriction of $\mathrm{\omega}$ to the curve
$F\left(  \sigma\right)  $ is a well-defined density. This density does not
change, if we replace $I$ by $-I.$ Therefore the Funk transform is well
defined. The function $M\mathrm{f}$ is also continuous. Suppose that the
condition (\textbf{iii}) is fulfilled. If a density $\mathrm{f}$ is supported
in a compact set $K\Subset X,$ then $M\mathrm{f}$ is supported in the compact
set $\Lambda\doteq\pi\left(  p^{-1}\left(  K\right)  \right)  \Subset\Sigma.$

\textbf{Example 1.} Let $X$ and $\Sigma$ be unit spheres in Euclidean 3-spaces
and the hypersurface $F$ be defined by the global incidence function $I\left(
x,\sigma\right)  =x_{1}\sigma_{1}+x_{2}\sigma_{2}+x_{3}\sigma_{3}.$ The
operator $M$ coincides with the classical Minkowski-Funk transform, \cite{F}.
The Funk transform can be also defined on projective planes $X=P^{2}%
=S^{2}/\mathbb{Z}_{2},\Sigma=P^{2}$, if we take $\pm I$ as local incidence
functions. The dual operator $M%
{{}^\circ}%
$ coincides with $M$ through the natural isomorphism $X\cong\Sigma.$

\textbf{Example 2.} Let $\left(  X,g\right)  $ be a Riemannian 2-manifold with
boundary and $\Sigma$ be the family of closed geodesic curves $\gamma.$ Take a
density $\mathrm{f}=f\mathrm{d}S$, where $\mathrm{d}S$ is the Riemannian area
form and $f$ is a continuous function with compact support. Then we can write
the geodesic integral transform as the Funk transform
\[
M\mathrm{f}\left(  \gamma\right)  =\int_{\gamma}f\mathrm{d}s,
\]
if we take an incidence function $I$ for the family $F$ such that $\left|
\nabla_{x}I\right|  =1.$ Any smooth weight function $w=w\left(  x,\sigma
\right)  $ can be included, by replacing the function $I$ to $w^{-1}I.$

\section{Above estimates}

The scale of Sobolev $L_{2}$-norms $\left\|  \cdot\right\|  ^{\alpha}%
,\alpha\in\mathbb{R}$ is defined for functions supported in an arbitrary
compact set $K\subset X.$ Fix a volume form $\mathrm{d}X$ in $X$ and define
$\left\|  \mathrm{f}\right\|  ^{\alpha}=\left\|  f_{0}\right\|  ^{ga}$ for a
density $\mathrm{f}=f_{0}\mathrm{d}X$ with support in $K.$ Denote by
$H_{K}^{\alpha}\left(  X,\Omega\right)  ,\;H_{K}^{\alpha}\left(  X\right)  $
the space of densities (distributions), respectively, of (generalized)
functions supported in $K$ with finite norm $\left\|  \cdot\right\|  _{\alpha
}$. For an arbitrary compact set $\Lambda\subset\Sigma$ we define the spaces
$H_{\Lambda}^{\alpha}\left(  \Sigma,\Omega\right)  ,\;H_{\Lambda}^{\alpha
}\left(  \Sigma\right)  $ in the same way.

\begin{proposition}
\label{p2}For any family $F$ that fulfils (\textbf{i-ii}), an arbitrary
compact set $K\Subset X$, any smooth function $\varepsilon$ in $\Sigma$ with
compact support and any real $\alpha$ the inequality holds
\[
\left\|  \varepsilon M\mathrm{f}\right\|  ^{\alpha+\left(  n-1\right)  /2}\leq
C_{\alpha}\left\|  \mathrm{f}\right\|  ^{\alpha}%
\]
for $\mathrm{f}\in H_{K}^{\alpha}\left(  X,\Omega\right)  $.
\end{proposition}

\textsf{Proof.} The Funk transform can be expressed as an oscillatory
integral
\[
M\mathrm{f}\left(  \sigma\right)  =\int_{K}\int_{\mathbb{R}}\exp\left(
2\pi\imath\tau I\left(  x,\sigma\right)  \right)  \mathrm{f}\left(  x\right)
\mathrm{d}\tau.
\]
The critical set of the phase function $\tau I\left(  x,\sigma\right)  $ is
the hypersurface $F\left(  \sigma\right)  $ and the condition $\mathrm{d}%
_{x}I\neq0$ implies that the phase function is non-degenerate. The
corresponding conic Lagrange manifold is
\[
\mathrm{L}=\{\left(  x,\sigma,\xi,\rho\right)  \in T^{\ast}\left(
X\times\Sigma\right)  ,I\left(  x,\sigma\right)  =0,\rho=\tau\mathrm{d}%
_{\sigma}I,\xi=\tau\mathrm{d}_{x}I,\tau\neq0\}.
\]

\begin{lemma}
Rank of the matrix $\frac{\partial\left(  x,\xi\right)  }{\partial\left(
\sigma,\rho\right)  }$ is equal to $2n$ in any point of $\Lambda.$
\end{lemma}

\textsf{Proof of Lemma.} Suppose that the rank is less $2n$. Then there exists
a vector $t=\left(  \delta x,\delta\sigma,\delta\xi,\delta\rho\right)  $ in
$T^{\ast}\left(  X\times\Sigma\right)  $ tangent to $\mathrm{L}$ such that
$\delta\sigma=0,\delta\rho=0.$ This yields
\begin{align*}
\mathrm{d}I\left(  \delta x\right)   &  =0,\delta\rho=\delta\tau
\mathrm{d}_{\sigma}I+\tau\mathrm{d}_{x}\mathrm{d}_{\sigma}I\left(  \delta
x\right)  ,\\
\delta\xi &  =\delta\tau\mathrm{d}_{x}I+\tau\mathrm{d}_{x}^{2}I\left(  \delta
x\right)
\end{align*}
for a tangent vector $\delta\tau$ to $\mathbb{R}.$ The first line implies that
the vector $\left(  \tau\delta x,\delta\tau\right)  $ fulfils the equation
$\left(  \tau\delta x,\delta\tau\right)  \Phi=0.$ By Proposition \ref{p1} this
vector vanishes, that is $\delta x=0,\delta\tau=0.$ The second line gives
$\delta\xi=0.$ $\blacktriangleright$

By this Lemma the projections of $\mathrm{L}$ to $T^{\ast}\left(  X\right)  $
and to $T^{\ast}\left(  \Sigma\right)  $ are submersions. In other terms,
$\mathrm{L}$ is locally the graph of a canonical transformation. The symbol
$a\left(  x,\sigma,\xi,\rho\right)  =1$ is a homogeneous function of $\xi
,\rho$ of order 0. The order $m$ of the Fourier integral operator $M$
satisfies the equation We have $m+\dim X\times\Sigma/4-N/2=0,$ where $\dim
X\times\Sigma=2n$ and $N=1$ is the number of variables $\tau.$ This yields
$m=\left(  1-n\right)  /2,$ which means that the functional
\[
I\left(  \psi\right)  \doteq\int_{\Sigma}\int_{X}\int_{\mathbb{R}}\exp\left(
2\pi\imath\tau I\left(  x,\sigma\right)  \right)  \psi\left(  x,\sigma\right)
\mathrm{d}\tau
\]
defined for a smooth densities $\psi$ in $X\times\Sigma$ with compact support,
is a distribution of the class $I^{\left(  1-n\right)  /2}\left(
X\times\Sigma,\mathrm{L}\right)  $ in the sense of Definition 25.4.9 of
\cite{H}. By Corollary 25.3.2 the operator $\varepsilon M$ defines a
continuous map $H_{K}^{\alpha}\left(  X,\Omega\right)  \rightarrow H_{\Lambda
}^{\alpha+\left(  n-1\right)  /2}\left(  \Sigma\right)  $ for any real
$\alpha,$ where $\Lambda=\mathrm{supp\,}e.$ $\blacktriangleright$

\begin{corollary}
\label{c1}If $F$ fulfils (\textbf{i}-\textbf{iii}), the Funk transform $M$ can
be extended to a bounded operator $H_{K}^{\alpha}\left(  X,\Omega\right)
\rightarrow H_{\Lambda}^{\alpha+\left(  n-1\right)  /2}\left(  \Sigma\right)
$ for any $\alpha,K$ and $\Lambda=\pi\left(  p^{-1}\left(  K\right)  \right)  .$
\end{corollary}

Due to (\textbf{iii}), we have\textbf{ }$\mathrm{supp\,}Mf\subset\Lambda$ and
the cutoff factor $\varepsilon$ can be dropped out.

\section{Dual Funk transform}

Let $\mathrm{\varphi}$ be a density in $\Sigma$ with compact support. Define
the dual Funk transform as follows
\begin{equation}
M%
{{}^\circ}%
\mathrm{\varphi}\left(  x\right)  =\lim_{\varepsilon\rightarrow0}\frac
{1}{2\varepsilon}\int_{\left|  I\left(  x,\cdot\right)  \right|
\leq\varepsilon}\mathrm{\varphi}=\int_{F\left(  x\right)  }\frac
{\mathrm{\varphi}}{\mathrm{d}_{\sigma}I}.\label{2}%
\end{equation}
If $F$ satisfies (\textbf{iii}), then the manifold $F\left(  x\right)  $ is
compact for any $x\in X$ and the integral (\ref{2}) is well defined for any
continuous density $\mathrm{\varphi}$ in $\Sigma.$ The natural pairing
\[
\left(  \mathrm{f},\mathrm{\phi}\right)  \mapsto\left\langle \mathrm{f}%
,\mathrm{\phi}\right\rangle \doteq\int_{X}\mathrm{f\phi}%
\]
is well defined for a densities $\mathrm{f}$\textrm{\ }and
functions\textrm{\ }$\mathrm{\phi}$ on $X$ provided one of them has compact support.

\begin{proposition}
\label{p5}The operator $-M%
{{}^\circ}%
$ is dual to $M.$
\end{proposition}

\textsf{Proof}\textsc{.} We have
\[
\left\langle M\mathrm{f},\mathrm{\varphi}\right\rangle =\int_{\Sigma}M\left(
\mathrm{f\,}\right)  \mathrm{\bar{\varphi}}=\int_{\Sigma}\mathrm{\bar{\varphi
}}\int_{F\left(  \sigma\right)  }\frac{\mathrm{f}}{\mathrm{d}_{x}I}=\int
_{F}\frac{\mathrm{f}\wedge\mathrm{\bar{\varphi}}}{\mathrm{d}_{x}I}=-\int
_{F}\mathrm{f}\wedge\frac{\mathrm{\bar{\varphi}}}{\mathrm{d}_{\sigma}I},
\]
since $\mathrm{d}I=\mathrm{d}_{x}I+\mathrm{d}_{\sigma}I=0$ on $F.$ The
right-hand side equals
\[
-\int_{X}\mathrm{f}\int_{F\left(  x\right)  }\frac{\mathrm{\bar{\varphi}}%
}{\mathrm{d}_{\sigma}I}=-\int_{X}\mathrm{f\,}M%
{{}^\circ}%
\left(  \mathrm{\bar{\varphi}}\right)  =-\left\langle \mathrm{f},M%
{{}^\circ}%
\mathrm{\varphi}\right\rangle .\blacktriangleright
\]

\section{Backprojection and two side estimates}

\textbf{Definition.} Fix some area forms $\mathrm{d}X$ in $X$ and
$\mathrm{d}\Sigma$ in $\Sigma$. The \textit{back projection} operator
\[
M^{\ast}:g\mapsto M%
{{}^\circ}%
\left(  g\mathrm{d}\Sigma\right)  \mathrm{d}X
\]
transforms functions defined in $\Sigma$ to densities in $X.$

\textbf{Definition. }We say that points $x,y\in X$ are conjugate with respect
to $F,$ if $x\neq y$ and the form $\mathrm{d}_{\sigma}I\left(  x,\sigma
\right)  \wedge\mathrm{d}_{\sigma}I\left(  y,\sigma\right)  $ defined in
$\Sigma$ vanishes.

\begin{theorem}
\label{t4}If a family $F$ has no conjugate points, fulfils (\textbf{i-ii}) and
the condition: \newline (\textbf{iv}) the projection $q:\pi^{-1}\left(
\Lambda\right)  \rightarrow G^{n-1}\left(  K\right)  $ is surjective for some
sets $K\Subset X,\Lambda\Subset\Sigma$, then for arbitrary cutoff function
$\varepsilon$ such that $\varepsilon=1$ in $\Lambda$ and any $\alpha>\beta$
the estimate
\begin{equation}
\left\|  \mathrm{f}\right\|  ^{\alpha}\leq C_{\alpha}\left\|  \varepsilon
M\mathrm{f}\right\|  ^{\alpha+\left(  n-1\right)  /2}+C_{\beta}\left\|
\mathrm{f}\right\|  ^{\beta}\label{9}%
\end{equation}
holds for the Funk transform of densities $\mathrm{f}$ supported in $K,$ where
$C_{\alpha}$ and $C_{\beta}$ do not depend on $f.$
\end{theorem}

\begin{lemma}
\label{l1}The composition $M^{\ast}\varepsilon M$ is an elliptic PDO in $K$ of
order $1-n$.
\end{lemma}

\textsf{Proof of Lemma.} Write $\mathrm{f}=f_{0}\mathrm{d}X,$ where $f_{0}$ is
a function supported in $K$ and calculate
\begin{align*}
\frac{M^{\ast}\varepsilon M\mathrm{f}}{\mathrm{d}X}\left(  y\right)   &
=\int_{F\left(  y\right)  }\frac{\varepsilon\left(  \sigma\right)
\mathrm{d}\Sigma}{\mathrm{d}_{\sigma}I}\int_{F\left(  \sigma\right)  }%
\frac{\mathrm{f}\left(  x\right)  }{\mathrm{d}_{y}I\left(  x,\sigma\right)
}\\
&  =\int_{I\left(  y,\sigma\right)  =0}\frac{\varepsilon\left(  \sigma\right)
\mathrm{d}\Sigma}{\mathrm{d}_{\sigma}I\left(  y,\sigma\right)  }\int_{I\left(
x,\sigma\right)  =0}\frac{\mathrm{d}X}{\mathrm{d}_{x}I\left(  x,\sigma\right)
}f_{0}\left(  x\right) \\
&  =-\int_{I\left(  y,\sigma\right)  =0}\int_{I\left(  x,\sigma\right)
=0}\frac{\varepsilon\left(  \sigma\right)  \mathrm{d}\Sigma}{\mathrm{d}%
_{\sigma}I\left(  y,\sigma\right)  \wedge\mathrm{d}_{\sigma}I\left(
x,\sigma\right)  }f_{0}\left(  x\right)  \mathrm{d}X,
\end{align*}
since $\mathrm{d}I=\mathrm{d}_{x}I+\mathrm{d}_{\sigma}I=0$ in $F.$ We can
write the right-hand side as $\int A\left(  y,x\right)  \mathrm{f}\left(
x\right)  ,$ where
\begin{equation}
A\left(  y,x\right)  =-\int_{F\left(  x\right)  \cap F\left(  y\right)  }%
\frac{\varepsilon\left(  \sigma\right)  \mathrm{d}\Sigma}{\mathrm{d}_{\sigma
}I\left(  y,\sigma\right)  \wedge\mathrm{d}_{\sigma}I\left(  x,\sigma\right)
}. \label{8}%
\end{equation}
The dominator does not vanishes for $y\neq x$, since there is no conjugate
points. The quotient $Q$ in (\ref{8}) is well defined as $n-2$-form up to an
additive term $\mathrm{d}_{\sigma}I\left(  y,\sigma\right)  \wedge
S+\mathrm{d}_{\sigma}I\left(  x,\sigma\right)  \wedge R$, where $S$ and $R$
are some $n-3$-forms. The integral of this term along the smooth manifold
$F\left(  x\right)  \cap F\left(  y\right)  $ vanishes and the function
$A\left(  y,x\right)  $ is a well defined and smooth, except for the diagonal.
Near the diagonal we can write $I\left(  y,\sigma\right)  =I\left(
x,\sigma\right)  +\sum\left(  y_{i}-x_{i}\right)  \partial I\left(
x,\sigma\right)  /\partial x_{i}+O\left(  \left|  y-x\right|  ^{2}\right)  $
and
\[
\mathrm{d}_{\sigma}I\left(  x,\sigma\right)  \wedge\mathrm{d}_{\sigma}I\left(
y,\sigma\right)  =\sum_{i}\left(  y_{i}-x_{i}\right)  \mathrm{d}_{\sigma
}I\left(  x,\sigma\right)  \wedge\mathrm{d}_{\sigma}\frac{\partial I\left(
x,\sigma\right)  }{\partial x_{i}}+O\left(  \left|  y-x\right|  ^{2}\right)
.
\]
The forms $\mathrm{d}_{\sigma}I\left(  x,\sigma\right)  \wedge\mathrm{d}%
_{\sigma}\partial I\left(  x,\sigma\right)  /\partial x_{i},i=1,...,n$ do not
vanish and are linearly independent, since of Proposition \ref{p1}. Therefore
the product $\mathrm{d}_{\sigma}I\left(  x,\sigma\right)  \wedge
\mathrm{d}_{\sigma}I\left(  y,\sigma\right)  $ is bounded by$\;c\left|
x-y\right|  $ from below as $y\rightarrow x.$ Therefore we have $A\left(
y,x\right)  =a\left(  y\right)  \left|  x-y\right|  ^{-1}+O\left(  1\right)
\;$near the diagonal, where $a$ is a smooth positive function. This implies
that $M^{\ast}\varepsilon M$ is a classical integral operator on $K$ with weak
singularity, moreover it is a pseudodifferential operator of order $1-n$. It
is an elliptic operator, since the symbol $a$ is positive.
$\blacktriangleright$

\textsf{Proof of Theorem.} The support of the function $M^{\ast}\varepsilon
M\mathrm{f}$ is contained in the compact set $p\left(  \mathrm{supp\,}%
\varepsilon\right)  \subset X.$ By Proposition \ref{p2} $M^{\ast}$ is $\left(
n-1\right)  /2$-smoothing operator, which yields
\begin{equation}
\left\|  M^{\ast}\varepsilon M\mathrm{f}\right\|  ^{\alpha+n-1}\leq C\left\|
\varepsilon M\mathrm{f}\right\|  ^{\alpha+\left(  n-1\right)  /2}. \label{6}%
\end{equation}
By Lemma \ref{l1} the operator $M^{\ast}\varepsilon M$ is elliptic of order
$n-1$, therefore the standard inequality holds
\[
\left\|  \mathrm{f}\right\|  ^{\alpha}\leq C_{\alpha}\left\|  M^{\ast
}\varepsilon M\mathrm{f}\right\|  ^{\alpha+n-1}+C_{\beta}\left\|
\mathrm{f}\right\|  ^{\beta}%
\]
for an arbitrary $\beta$ and some constants $C_{\alpha},C_{\beta}.$ Taking in
account (\ref{6}) yields (\ref{9}). $\blacktriangleright$

\begin{corollary}
The eigenvalues $\lambda_{k}$ of the operator $M^{\ast}\varepsilon M$ numbered
in decreasing order satisfy the estimate
\[
ck^{\left(  1-n\right)  /2}\leq\lambda_{k}\leq Ck^{\left(  1-n\right)
/2},\,k\geq k_{0}.
\]
\end{corollary}

For the Radon transform the eigenvalues are calculated in \cite{N}.

\begin{corollary}
\label{c3}Suppose that for some $\beta<\alpha$ the equation $M\mathrm{f}%
=0,\mathrm{f}\in H_{K}^{\beta}\left(  X,\Omega\right)  $ implies
$\mathrm{f}=0.$ Then the two-side estimate holds:
\begin{equation}
c_{\alpha}\left\|  \mathrm{f}\right\|  ^{\alpha}\leq\left\|  \varepsilon
M\mathrm{f}\right\|  ^{\alpha+\left(  n-1\right)  /2}\leq C_{\alpha}\left\|
\mathrm{f}\right\|  ^{\alpha}.\label{7}%
\end{equation}
\end{corollary}

\textsf{Proof. }The right-hand side inequality follows from Proposition
\ref{p2}. Suppose that the left-hand side estimate does hold for no
$c_{\alpha}.$ Then there exists a sequence $\left\{  \mathrm{f}_{k}\right\}
\subset H_{K}^{\alpha}\left(  X,\Omega\right)  $ such that
\begin{equation}
\left\|  \mathrm{f}_{k}\right\|  _{\alpha}\geq k\left\|  M\mathrm{f}%
_{N}\right\|  _{\alpha+1/2},\left\|  \mathrm{f}_{k}\right\|  _{\beta
}=1,k=1,2,... \label{4}%
\end{equation}
The inequality (\ref{9}) implies that $\left\|  \mathrm{f}_{k}\right\|
_{\alpha}\leq2C_{\beta}$ for $k>2C_{\alpha}$ and $\left\|  M\mathrm{f}%
_{k}\right\|  _{\alpha+1/2}\rightarrow0.$ Because the imbedding $H_{K}%
^{\alpha}\left(  X,\Omega\right)  \rightarrow H_{K}^{\beta}\left(
X,\Omega\right)  $ is compact, we can choose a subsequence (denote it again
$\left\{  \mathrm{f}_{k}\right\}  $) such that $\mathrm{f}_{k}\rightarrow
\mathrm{g}$ in $H_{K}^{\beta}\left(  X,\Omega\right)  .$ By Proposition
\ref{p2} $\left\|  M\mathrm{f}_{k}\rightarrow M\mathrm{g}\right\|
_{\beta+\left(  n-1\right)  /2}\rightarrow0,$ which implies $M\mathrm{g}=0.$
By the condition $\mathrm{g}=0;$ it follows that $\left\|  \mathrm{f}%
_{k}\right\|  _{\beta}\rightarrow0$ in contradiction with (\ref{4}).
$\blacktriangleright$

\textbf{Remarks. }Mukhometov's result \cite{M} implies the estimate $\left\|
\mathrm{f}\right\|  ^{0}\leq C\left\|  M\mathrm{f}\right\|  ^{1}$ for the case
$n=2$. An estimate of this kind for more general situation was obtained by
Sharafutdiniv \cite{Sh}, Ch. IV. Inequalities for Sobolev norms are well known
for the Radon transform. Estimates for shift derivatives of order $\alpha+1/2$
$(n=2)$ were obtained by several authors. Natterer \cite{N} has shown that
(\ref{7}) holds also for angular derivatives. For the attenuated Radon
transform see Rullgard \cite{R}.

Our approach is similar to that of Lavrent'ev and Bukhgeim \cite{LB}, where
the composition $M^{\ast}M$ was described as an integral operator in the local
case. Guillemin \cite{G},\cite{GS} has defined the `generalized Radon
transform' $R$ for an arbitrary double fibration. This transfrom is treated as
an elliptic Fourier integral operator and $R^{\ast}R$ is shown to be a
pseudodifferential elliptic operator under the `Bolker condition'. This
condition is equivalent to absence of conjugate points in our situation. More
details are given in the paper of T. Quinto \cite{Q}.

\section{Range conditions and approximation}

Let $K$ be a compact set in $X$ and $\alpha\in\mathbb{R.}$ We define
$H^{\alpha}\left(  K,\Omega\right)  \,$to be the dual space of $H_{K}%
^{-\alpha}\left(  X\right)  $ and use the notation $\left\|  \cdot\right\|
^{\alpha}$ for the norm in $H^{\alpha}\left(  K,\Omega\right)  .$ The trace
operator $H_{L}^{\alpha}\left(  X,\Omega\right)  \rightarrow H^{\alpha}\left(
K,\Omega\right)  $ is well defined and is open for an arbitrary compact set
$L\subset X$ such that $K\Subset L,$ since $H_{K}^{-\alpha}\left(  X\right)  $
is a subspace of $H_{L}^{-\alpha}\left(  X\right)  .$ Therefore $H^{\alpha
}\left(  K,\Omega\right)  $ can be realized as the quotient space of
$H_{L}^{\alpha}\left(  X,\Omega\right)  $ modulo the kernel of the trace
operator. The last one consists of densities $f$ supported in $L\backslash K.$
For any $\beta>\alpha$ we have the operator $\eta^{\prime}:H^{\beta}\left(
K,\Omega\right)  \rightarrow H^{\alpha}\left(  K,\Omega\right)  ,$ which is
dual to the natural imbedding $\eta:H_{K}^{-\alpha}\left(  X\right)
\rightarrow H_{K}^{-\beta}\left(  X\right)  .$ If the boundary of $K$ is
smooth, the imbedding $\eta$ has dense image and $\eta^{\prime}$ is injective.
Then we can define the intersection $H^{\infty}\left(  K,\Omega\right)
\doteq\cap_{\alpha}H^{\alpha}\left(  K,\Omega\right)  ;$ any density
$\mathrm{f}\in H^{\infty}\left(  K,\Omega\right)  $ is smooth in the interior
of $K$. Similarly, we define $H^{\alpha}\left(  K\right)  \doteq\left(
H_{K}^{-\alpha}\left(  X,\Omega\right)  \right)  ^{\prime}.$

Suppose that the incidence manifold $F$ fulfils (\textbf{i-iii}). If a density
$\mathrm{f}$ is supported in a compact set $K,$ then the support of
$M\mathrm{f}$ is contained in the compact set $\Lambda=\pi\left(
p^{-1}\left(  K\right)  \right)  .$ The hypersurface $F\left(  x\right)  $ is
compact for any point $x\in X$ and the dual transform $M%
{{}^\circ}%
$ is well defined for all continuous densities in $\Sigma.$ Moreover, it can
be extended to a continuous  operator $M%
{{}^\circ}%
:H^{\alpha}\left(  \Lambda,\Omega\right)  \rightarrow H^{\alpha+\left(
n-1\right)  /2}\left(  K\right)  $ for any $\alpha$ by means of the duality
\[
\left\langle M%
{{}^\circ}%
\mathrm{g},\mathrm{f}\right\rangle =-\left\langle \mathrm{g},M\mathrm{f}%
\right\rangle ,\;\mathrm{f}\in H_{K}^{-\alpha-\left(  n-1\right)  /2}\left(
X,\Omega\right)  ,\,\mathrm{g}\in H^{\alpha}\left(  \Lambda,\Omega\right)  .
\]
By Proposition \ref{p2}, we have $M\mathrm{f}\in H_{\Lambda}^{-\alpha}\left(
\Sigma\right)  ,$ hence the right-hand side is well defined. If $\alpha\geq0,$
the function $M%
{{}^\circ}%
\mathrm{g}$ defined by this formula is equal to the integral (\ref{2}), which
has sense, at least, for almost all $x\in K.$

\begin{theorem}
\label{t3}Suppose that $F$ satisfies (\textbf{i}-\textbf{iii}) and has no
conjugate points. Then for any $K\Subset X$ and arbitrary $\alpha
\in\mathbb{R\cup}\left\{  \infty\right\}  $ the image of the Funk operator
\[
M:H_{K}^{\alpha}\left(  X,\Omega\right)  \rightarrow H_{\Lambda}%
^{\alpha+\left(  n-1\right)  /2}\left(  \Sigma\right)  ,\;\Lambda=\pi\left(
p^{-1}\left(  K\right)  \right)
\]
is closed and coincides with the subspace of functions $\varphi\in H_{\Lambda
}^{\alpha+\left(  n-1\right)  /2}\left(  \Sigma\right)  $ such that
$\int\mathrm{g}\varphi=0$ for any solution $\mathrm{g}\in H^{-\alpha-\left(
n-1\right)  /2}\left(  \Lambda,\Omega\right)  $ of the equation
\begin{equation}
M%
{{}^\circ}%
\mathrm{g}\left(  x\right)  =0,x\in K.\label{5}%
\end{equation}
\end{theorem}

\textsf{Proof.} The image of $\mathrm{\,}M$ is closed by Theorem \ref{t4},
thereby it coincides with the polar of the kernel of the dual operator $M%
{{}^\circ}%
.$ $\blacktriangleright$

\begin{theorem}
\label{t5}If $F$ fulfils (\textbf{i}-\textbf{iii}) and has no conjugate
points. Then for any set $K\Subset X$ with smooth boundary and arbitrary real
$\alpha\in\mathbb{R}$ any density $\mathrm{g}\in H^{\alpha}\left(
\Lambda,\Omega\right)  ,\Lambda=\pi\left(  p^{-1}\left(  K\right)  \right)  $
that fulfils (\ref{5}) can be approximated by solutions $\mathrm{h}\in
H^{\infty}\left(  \Lambda,\Omega\right)  \doteq\cap_{\alpha}H^{\alpha}\left(
\Lambda,\Omega\right)  .$
\end{theorem}

\textsf{Proof.} Let $\mathrm{Sol}^{\beta}$ denote the space of solutions of
(\ref{5}) in the class $H^{\beta}\left(  \Lambda,\Omega\right)  .$ We show
first that $\mathrm{g}$ can be approximated by elements of $\mathrm{Sol}%
^{\beta}$ for any $\beta>\alpha.$ It is sufficient to check that any
functional $\mathrm{\phi}$ on $H^{\alpha}\left(  \Lambda,\Omega\right)  $ that
is equal to zero on $\mathrm{Sol}^{\beta}$ also vanishes on $\mathrm{g.}$ The
dual space is isomorphic to $H_{\Lambda}^{-\alpha}\left(  \Sigma\right)  ,$
which implies $\mathrm{\phi}\in$ $H_{\Lambda}^{-\alpha}\left(  \Sigma\right)
.$ By Corollary \ref{c1} the Funk transform defines the continuous operator
\[
M_{\beta}:H_{K}^{-\beta-\left(  n-1\right)  /2}\left(  X,\Omega\right)
\rightarrow H_{\Lambda}^{-\beta}\left(  \Sigma\right)  .
\]
By Theorem \ref{t3}, the range of this operator is closed and coincides with
the polar set of $\mathrm{Sol}^{\beta}$. It follows that, $\mathrm{\phi
}=M\mathrm{\psi}$ for a density $\mathrm{\psi}\in H_{K}^{-\beta-\left(
n-1\right)  /2}\left(  X,\Omega\right)  .$ By Theorem \ref{t4} we have
$\mathrm{\psi}\in H_{K}^{-\alpha-\left(  n-1\right)  /2}\left(  X,\Omega
\right)  $ and can write
\[
\left\langle \mathrm{\phi},\mathrm{g}\right\rangle =\left\langle
M\mathrm{\psi},\mathrm{g}\right\rangle =-\left\langle \mathrm{\psi},M%
{{}^\circ}%
\mathrm{g}\right\rangle =0,
\]
since the function $M%
{{}^\circ}%
\mathrm{g}$ is well defined as element of $H_{K}^{\alpha+\left(  n-1\right)
/2}\left(  X\right)  .$ This yields that $\mathrm{g}$ is contained in the
closure of the space $\mathrm{Sol}^{\beta}.$

Now we approximate $\mathrm{g}$ by elements of the space $H^{\infty}\left(
\Lambda,\Omega\right)  .$ Let $\left\|  \cdot\right\|  _{\Lambda}^{\alpha}$ be
the norm in the space $H^{\alpha}\left(  \Lambda,\Omega\right)  ,$ which is
dual to the norm $\left\|  \cdot\right\|  ^{-\alpha}.$ We may assume that the
norm $\left\|  \cdot\right\|  _{\Lambda}^{\alpha}$ is monotone increasing
function of $\alpha.$ Take an arbitrary $\varepsilon$ and choose a function
$\mathrm{h}_{1}\in H^{\alpha+1}\left(  \Lambda,\Omega\right)  $ such $\left\|
\mathrm{h}_{1}-\mathrm{g}\right\|  _{\Lambda}^{\alpha}<\varepsilon/2,$ then we
choose a function $\mathrm{h}_{2}\in H^{\alpha+2}\left(  \Lambda
,\Omega\right)  $ such that $\left\|  \mathrm{h}_{2}-\mathrm{h}_{1}\right\|
_{\Lambda}^{\alpha+1}<\varepsilon/4$ and so on. We obtain a sequence $\left\{
\mathrm{h}_{k}\right\}  $ such that $\left\|  \mathrm{h}_{k+1}-\mathrm{h}%
_{k}\right\|  _{\Lambda}^{\alpha+k}<2^{-k-1}\varepsilon$ for $k=1,2,...$ This
sequence converges to an element $\mathrm{h}$ in any space $H^{\beta}\left(
\Lambda,\Omega\right)  ,\beta>\alpha.$ It follows that $\mathrm{h}\in
H^{\infty}\left(  \Lambda,\Omega\right)  $ and $\left\|  \mathrm{h}%
-\mathrm{g}\right\|  _{\Lambda}^{\alpha}\leq\varepsilon.\blacktriangleright$

\begin{corollary}
Under conditions of Theorem \ref{t3}, it is sufficient to check the equation
$\int\mathrm{g}\varphi=0$ for densities $\mathrm{g}\in H^{\infty}\left(
\Lambda,\Omega\right)  $ that satisfies (\ref{5}).
\end{corollary}

\section{Thermoacoustic tomography}

We apply the above results for the thermo/opto/photoacoustic geometry. First,
consider the case of complete acquisition geometry. Let $X$ be the open unit
ball in an Euclidean space $\mathbf{E},$ $S_{R}$ be the sphere of radius $R>1$
and $\Sigma=S_{R}\times\mathbb{R}$. The manifold $F\subset X\times\Sigma$ is
given by the equation $I\left(  x;y,r\right)  =\left|  y-x\right|  -r=0,y\in
S_{R},0<r.$ The manifold $F$ obviously fulfils (\textbf{i),(iii}) and has no
conjugate points. Check that the condition (\textbf{ii}) is also satisfied. It
is sufficient to check that, the sphere $F\left(  y,r\right)  $ can not be
tangent to $F\left(  z,s\right)  $ at a point $x\in X,$ if the points $\left(
y,r\right)  $ and $\left(  z,s\right)  $ are sufficiently close in $\Sigma.$
The condition $\left|  y+z\right|  >2$ is sufficient for this. The Funk
operator
\begin{equation}
M\mathrm{f}\left(  y,r\right)  =\int_{\left|  x-y\right|  =r}f_{0}\left(
x\right)  \mathrm{d}S,\;\mathrm{f}=f_{0}\mathrm{d}x \label{mfs}%
\end{equation}
is the spherical integral transform, where $\mathrm{d}S$ is the Euclidean
surface area form on spheres. The kernel of dual transform $M%
{{}^\circ}%
$ consists of densities $\mathrm{\varphi}=\phi\mathrm{d}S\mathrm{d}r$ in
$\Sigma$ such that
\[
0=\int_{F\left(  x\right)  }\frac{\mathrm{\varphi}}{\mathrm{d}_{\sigma}I}%
=\int\frac{\phi\mathrm{d}S\mathrm{d}r}{\mathrm{d}\left(  \left|  x-y\right|
-r\right)  }=-\int_{F\left(  x\right)  }\phi\mathrm{d}S.
\]
for any $x\in X.$ Theorem \ref{t5} yields

\begin{corollary}
\label{c4}For any compact set $K\subset X$ with smooth boundary and arbitrary
$\alpha\in\mathbb{R}$ $\cup\left\{  \infty\right\}  $ the range of the Funk
operator $M:H_{K}^{\alpha}\left(  X,\Omega\right)  \rightarrow H_{\Lambda
}^{\alpha+\left(  n-1\right)  /2}\left(  \Sigma\right)  $ coincides with the
set of functions $g$ in $\Sigma$ such that
\begin{equation}
\int_{\Sigma}g\mathrm{\varphi}=0\label{gph}%
\end{equation}
for any density $\mathrm{\varphi}=\phi\mathrm{d}S\mathrm{d}r$ such that
$\phi\in C^{\infty}\left(  \Sigma\right)  $ and
\begin{equation}
\int_{F\left(  x\right)  }\phi\mathrm{d}S=0,\;x\in K.\label{ifs}%
\end{equation}
\end{corollary}

\textbf{Remark. }For the operator $M$ acting on $C^{\infty}$-densities the
range conditions were given in the papers \cite{Pa}, \cite{AKQ}, \cite{Fi}.
The conditions of \cite{AKQ} and \cite{Fi} give full description of the range
of $M,$ but have implicit form.

We extract some explicit range conditions from Corollary \ref{c4}. For an
arbitrary $x\in K,$ the manifold $F\left(  x\right)  $ is the intersection of
the cone surface $\left|  y-x\right|  =r$ with the cylinder $\Sigma.$ This
intersection is contained in the hyperplane
\[
P\left(  x\right)  =\left\{  y,s;2\left\langle x,y\right\rangle +s=\left|
x\right|  ^{2}+R^{2}\right\}  \subset\mathbf{E}\times\mathbb{R,}%
\]
where we set $s=r^{2}.$ Thus, the condition (\ref{ifs}) means vanishing of
integrals of $\phi\mathrm{d}S$ over intersections of $\Sigma$ with the
hyperplanes $P\left(  x\right)  ,x\in K\mathbb{.}$ Suppose that $\phi$ is a
polynomial in $s:\phi\left(  y,s\right)  =\sum\phi_{k}\left(  y\right)
\left(  s-R^{2}\right)  ^{k}$ and have%
\[
\int_{F\left(  x\right)  }\phi\mathrm{d}S=\sum_{k}\int_{S}\phi_{k}\left(
y\right)  \left(  \left|  x\right|  ^{2}-2\left\langle x,y\right\rangle
\right)  ^{k}\mathrm{d}S\left(  y\right)  .
\]
Set $x=tz$ for $\left|  z\right|  =1$ and $0\leq t<1$ and develop the
right-hand side in powers of $t:$%

\begin{align*}
\int_{F\left(  x\right)  }\phi\mathrm{d}S  &  =\sum\int_{S}\phi_{k}\left(
y\right)  \left(  t^{2}-2t\left\langle z,y\right\rangle \right)
^{k}\mathrm{d}S\\
&  =\int\phi_{0}\left(  y\right)  \mathrm{d}S-t\int2\phi_{1}\left(  y\right)
\left\langle z,y\right\rangle \mathrm{d}S\\
&  +t^{2}\int\left[  \phi_{1}\left(  y\right)  +4\phi_{2}\left(  y\right)
\left\langle z,y\right\rangle ^{2}\right]  \mathrm{d}S\\
&  -t^{3}\int\left[  4\phi_{2}\left(  y\right)  \left\langle z,y\right\rangle
+8\phi_{3}\left(  y\right)  \left\langle z,y\right\rangle ^{3}\right]
\mathrm{d}S\\
&  +t^{4}\int\left[  \phi_{2}\left(  y\right)  +12\phi_{3}\left(  y\right)
\left\langle z,y\right\rangle ^{2}+16\phi_{4}\left(  y\right)  \left\langle
z,y\right\rangle ^{4}\right]  \mathrm{d}S\\
+..  &  =0
\end{align*}
The right-hand side vanishes for all $t,$ which yields the system of
equations
\begin{align*}
\int\phi_{0}\left(  y\right)  \mathrm{d}S  &  =0,\\
\int\phi_{1}\left(  y\right)  \left\langle z,y\right\rangle \mathrm{d}S  &
=0,\\
\int\left[  4\phi_{2}\left(  y\right)  \left\langle z,y\right\rangle ^{2}%
+\phi_{1}\left(  y\right)  \right]  \mathrm{d}S  &  =0,\\
\int\left[  2\phi_{3}\left(  y\right)  \left\langle z,y\right\rangle ^{3}%
+\phi_{2}\left(  y\right)  \left\langle z,y\right\rangle \right]  \mathrm{d}S
&  =0,\\
\int\left[  16\phi_{4}\left(  y\right)  \left\langle z,y\right\rangle
^{4}+12\phi_{3}\left(  y\right)  \left\langle z,y\right\rangle ^{2}+\phi
_{2}\left(  y\right)  \right]  \mathrm{d}S  &  =0,\\
&  ...
\end{align*}

\begin{corollary}
Any solution $\left(  \phi_{0},\phi_{1},\phi_{2},...\right)  $ of this system
such that $\phi_{j}=0$ for all $j>k$ for some $k$ yields a function
$\phi\left(  y,s\right)  =\sum_{0}^{k}\phi_{j}\left(  y\right)  \left(
s-R^{2}\right)  ^{j}\;$that is a polynomial in $s$ of order $k$, fulfils
(\ref{ifs}) and is orthogonal to the range of $M.$
\end{corollary}

There are many solutions of this form, since the system has triangle form with
diagonal terms
\[
\int\phi_{j}\left(  y\right)  \left\langle z,y\right\rangle ^{j}%
\mathrm{d}S,\,\left|  z\right|  =1,j=0,1,...,k.
\]
To solve these equation we only need to fix the moments of $\phi_{j}$ of
degree $j$. There are only $%
\genfrac{(}{)}{}{}{n+j-1}{n-1}%
$ linearly independent $j$-moments, hence one can find infinitely many
independent solutions which are finite sums of harmonics. In particular, we
can take for $\phi_{0}$ any function on the sphere with zero average, an
arbitrary function $\phi_{1}$ with zero linear moments of $\phi_{1}$ and set
$\phi_{k}=0$ for $k>1$ etc. The range conditions of S. Patch \cite{Pa} are
apparently contained in (\ref{gph}) for polynomial $\phi.$

\section{Partial scan and Kaczmarz method}

In the case of the partial scan geometry the analysis is more complicated. Let
again $X$ be the open ball of radius $1$ in $\mathbf{E}$ and $\Sigma_{\delta
}=S_{\delta}\times\mathbb{R}_{+}$, where $S_{\delta}\doteq$ $\left\{
y;\left|  y\right|  =R,y_{1}>-\delta\right\}  $ for some $\delta>0.$ The
manifold $F_{\delta}$ is defined in $X\times\Sigma_{\delta}$ by the same
incidence function $I,$ that is, $F_{\delta}\left(  y,r\right)  $ is the
sphere of radius $r$ with the center $y\in S_{\delta}.$ The manifold
$F_{\delta}$ fulfils (\textbf{i}),(\textbf{ii)} and has no conjugate points,
but does not fulfil the condition (\textbf{iii}). On the other hand,
$F_{\delta}$ satisfies (\textbf{iv}) for the unit half-ball $K=\left\{
x,\left|  x\right|  \leq1,x_{1}\geq0\right\}  $ and $\Lambda=\left\{
y,r:R-1\leq r\leq R+1\right\}  $. Consider the Funk transform
\begin{equation}
M:H_{K}^{\alpha}\left(  X,\Omega\right)  \rightarrow H_{\Lambda}%
^{\alpha+\left(  n-1\right)  /\alpha}\left(  \Sigma\right)  \label{mbl}%
\end{equation}
defined as in (\ref{mfs}). Take a smooth function $\varepsilon_{0}\geq0$ on
$\mathbb{R}$ supported in the interval $(-\delta,\infty)$ such that
$\varepsilon_{0}\left(  t\right)  =1$ for $0\leq t\leq1.$ Set $\varepsilon
\left(  y\right)  =\varepsilon_{0}\left(  y_{1}\right)  $ and consider the
operator $M^{\ast}\varepsilon M.$

\begin{proposition}
\label{p6}The operator$\mathrm{\,}$(\ref{mbl}) is injective for any $\delta>0$
and arbitrary $\alpha>1/2.$ The inequality (\ref{7}) holds for the operator
$M^{\ast}\varepsilon M$ and any $\alpha>1/2.$
\end{proposition}

$\blacktriangleleft$ We prove that the equation $M\mathrm{f}=0$ in
$\Sigma_{\varepsilon}$ for a density $\mathrm{f}\in H_{K}^{\alpha}\left(
X,\Omega\right)  $ implies $\mathrm{f}\neq0.$ This condition means that the
spherical means of $\mathrm{f}$ vanish for spheres centered at points $y\in
S_{\delta}.$ By the Lin-Pinkus theorem \cite{LP} this implies that either
$\mathrm{f}=0$ or a non-trivial harmonic polynomial $h$ vanishes on
$S_{\delta}$ (the continuity condition for $\mathrm{f}$ in \cite{LP} can be
weakened). The last option is impossible, since $S_{\delta}$ is strictly
convex and $\Delta h$ can not vanish identically near $S_{\delta}.$ The
estimate (\ref{7}) now follows from Theorem \ref{t4} and Corollary \ref{c3}.
$\blacktriangleright$

We show that the Kaczmarz method can be adapted for inversion of the operator
(\ref{mbl}) as well as for a general Funk operator $M$. Let $K,\Lambda$ be
compact manifolds, occasionally with boundaries, $F$ be a closed hypersurface
in $K\times\Lambda$ that fulfils (\textbf{i}) and (\textbf{ii}). We want to
find a solution $\mathrm{f}\in H^{0}\left(  X,\Omega\right)  $ of the
equation
\begin{equation}
M\mathrm{f}=\mathrm{\varphi}\label{12}%
\end{equation}
for a function $\mathrm{\varphi}\in H^{0}(\Lambda).$ By Proposition \ref{p2},
the left-hand side is contained in the space $H^{\left(  n-1\right)
/2}\left(  \Lambda\right)  .\subset H^{0}\left(  \Lambda\right)  .$

\textbf{Example 3.} Take for $K$ the closed unit ball in $\mathbf{E}$, for
$\Lambda$ the product $\bar{S}_{\delta}\times\left[  R-1,R+1\right]  $ and for
$F$ the manifold of spheres as above.\textbf{ }The conditions (\textbf{i}0 and
(\textbf{ii}) are fulfilled.

Fix a volume form $\mathrm{d}X$ in $K,$ a volume form $\mathrm{d}\Sigma$ in
$\Lambda$ and consider the operator $MM^{\ast}:H^{0}\left(  \Lambda\right)
\rightarrow H^{0}\left(  \Lambda\right)  .$ It is non-positive; set
$R=-MM^{\ast}+\theta I,$ where $I$ is the identity operator and $\theta>0.$
The operator $R$ is self-adjoint, positive and invertible. Following \cite{N},
we choose a real parameter $\omega$ and set $Q\doteq I-\omega M^{\ast}%
R^{-1}M.$ We use the notation $\left\|  \cdot\right\|  =\left\|
\cdot\right\|  ^{0}.$

\begin{lemma}
\label{l2}We have $\left\|  Q\mathrm{g}\right\|  <\left\|  \mathrm{g}\right\|
$ for $0<\omega<2$ and any $\mathrm{g}\in H^{0}\left(  K,\Omega\right)  $ such
that $M\mathrm{g}\neq0$
\end{lemma}

\textsf{Proof.} We have by Proposition \ref{p5}
\begin{align*}
\left\|  Q\mathrm{g}\right\|  ^{2}  &  =\left\|  \mathrm{g}\right\|
^{2}-2\omega\left\langle \mathrm{g},M^{\ast}R^{-1}M\mathrm{g}\right\rangle
+\omega^{2}\left\langle M^{\ast}R^{-1}M\mathrm{g},M^{\ast}R^{-1}%
M\mathrm{g}\right\rangle \\
&  =\left\|  \mathrm{g}\right\|  ^{2}-2\omega\left\langle M\mathrm{g}%
,R^{-1}M\mathrm{g}\right\rangle -\omega^{2}\left\langle R^{-1}M\mathrm{g}%
,M^{\ast}MR^{-1}M\mathrm{g}\right\rangle \\
&  =\left\|  \mathrm{g}\right\|  ^{2}-2\omega\left\langle M\mathrm{g}%
,R^{-1}M\mathrm{g}\right\rangle +\omega^{2}\left\langle R^{-1}M\mathrm{g}%
,M\mathrm{g}\right\rangle -\omega^{2}\left\langle R^{-1}M\mathrm{g}%
,R^{-1}M\mathrm{g}\right\rangle \\
&  \leq\left\|  \mathrm{g}\right\|  ^{2}-\omega\left(  2-\omega\right)
\left\langle M\mathrm{g},R^{-1}M\mathrm{g}\right\rangle ,
\end{align*}
since $\left\langle R^{-1}M\mathrm{g},R^{-1}M\mathrm{g}\right\rangle \geq0.$
The term $\left\langle M\mathrm{g},R^{-1}M\mathrm{g}\right\rangle $ is
positive, if $M\mathrm{g}\neq0.$ $\blacktriangleright$

Take an arbitrary density $\mathrm{f}^{0}$ and construct the sequence
$\mathrm{f}^{k},k=1,2,...$ by means of the recurrent formula
\[
\mathrm{f}^{k+1}=\mathrm{f}^{k}+\omega M^{\ast}R^{-1}\left(  \mathrm{\varphi
}-M\mathrm{f}^{k}\right)  .
\]

\begin{theorem}
If $M$ is injective and $\mathrm{\varphi}$ fulfils the range conditions, we
have $\mathrm{f}^{k}\rightarrow\mathrm{f}$, where $\mathrm{f}$ is a solution
of (\ref{12}).
\end{theorem}

\textsf{Proof.} We have
\begin{align*}
Q\left(  \mathrm{f}^{k}-\mathrm{f}\right)   &  =\mathrm{f}^{k}-\omega M^{\ast
}R^{-1}M\mathrm{f}^{k}-\mathrm{f}+\omega M^{\ast}R^{-1}M\mathrm{f}\\
&  =\mathrm{f}^{k}+\omega M^{\ast}R^{-1}\left(  \mathrm{\varphi}%
-M\mathrm{f}^{k}\right)  -\mathrm{f}=\mathrm{f}^{k+1}-\mathrm{f.}%
\end{align*}
It follows that
\[
\left\|  \mathrm{f}^{k+1}-\mathrm{f}\right\|  <\left\|  \mathrm{f}%
^{k}-\mathrm{f}\right\|  <...<\left\|  \mathrm{f}^{0}-\mathrm{f}\right\|
\]
and $\mathrm{f}^{k}\rightarrow\mathrm{g}$ strongly in $L_{2}\left(
K,\Omega\right)  .$ We have $\left\|  Q\left(  \mathrm{g}-\mathrm{f}\right)
\right\|  =\left\|  \mathrm{g}-\mathrm{f}\right\|  ,$ which yields
$\mathrm{g}=\mathrm{f}$ by Lemma \ref{l2}. $\blacktriangleright$

Another inversion method is developed by Popov and Sushko \cite{PS}.

\end{document}